\documentclass{article}
\usepackage{amssymb,amsmath,amsfonts}
\usepackage{epsfig}
\usepackage{enumerate}
\newtheorem{proposition}{Proposition}
\newtheorem{theorem}{Theorem}
\begin{document}
\title{Boundary integral solution of
potential problems arising in the modelling of electrified oil
films}
\author{David J.\ Chappell\\
School of Science and Technology,\\
Nottingham Trent University,\\
Clifton Campus,\\
Clifton Lane,\\
Nottingham, UK\\
NG11 8NS\\
david.chappell@ntu.ac.uk} \maketitle
\begin{abstract}
We consider a class of potential problems on a periodic half-space
for the modelling of electrified oil films, which are used in the
development of novel switchable liquid optical devices (diffraction
gratings). A boundary integral formulation which reduces the problem
to the study of the oil-air interface alone is derived and solved in
a highly efficient manner using the Nystr\"{o}m method. The oil
films encountered experimentally are typically very thin and thus an
interface-only integral representation is important for avoiding the
near-singularity problems associated with boundary integral methods
for long slender domains.  The super-algebraic convergence of the
proposed methods is discussed and demonstrated via appropriate
numerical experiments.
\end{abstract}

\section{Introduction}

We consider a transmission problem for the Laplace equation in a
periodic half-plane for modelling the electric potential on an
interface between a thin layer of oil and the surrounding air.  A
fast solution algorithm is sought for the case where the
time-dependent interface position is given as the (discrete) data
generated from the solution of an associated coupled thin fluid and
interface dynamics problem. In particular, the interface motion is a
result of dielectrophoresis forces, where the potential is applied
via electrodes placed at the base of the oil layer.
Dielectrophoresis is defined as the motion of matter caused by
polarization effects in a nonuniform electric field
\cite{HP58,HP78}. The use of dielectrophoresis forces to stimulate
fluid motion has been considered in \cite{TJ01}, and is of interest
in optics due to potential applications in the development of novel
switchable liquid optical devices (diffraction gratings)
\cite{CB09}. In this work we focus on the solution of the
transmission problem for the Laplace equation, and solve for the
electric potential on the oil-air interface. We also compute the
normal and tangential derivatives of the potential on the interface
since these are the quantities of interest for the boundary data of
the associated thin film fluid problem.

The work here aims to provide a first step in moving from the study
of the static problem described in \cite{CB11} to the dynamic
situation described above, via an efficient solution method for the
underlying potential problem. Related work on the simulation of
electrified fluids has been carried out using a range of techniques
including asymptotic approaches \cite{CB11,LY06}, level set methods
\cite{SW06}, finite element approaches for coupled fluid flow and
dynamic interface models \cite{SW13}, and boundary integral methods
for coupled potential and dynamic interface problems \cite{DV05}.

Our strategy will be to reformulate our transmission problem in a
periodic half-plane as a boundary integral problem defined only on
the interface. A commonly cited advantage of boundary integral
methods is the reduction in dimensionality (by one) to the boundary
of the domain being studied. Here, by making explicit use of the
half-space solution of a related boundary value problem we go one
step further and reduce our study to the interface part of the
boundary. Considering that our goal is the computation of the
potential and its normal and tangential derivatives on the
interface, this reduction of the problem provides the ideal platform
for a highly efficient method. A major advantage over a conventional
boundary integral formulation where the entire boundary of the
finite domain is discretized, is that the near-singularity problems
due to the long slender geometry are completely avoided here.

Two further major reasons for
adopting a boundary integral approach are that the relevant Green's
functions are available in a simple closed form making it relatively
simple to implement a boundary integral method (compared to say
\cite{DV05}, where the periodic Green's function must be
approximated via fast summation methods). Secondly, the infinite
domains are dealt with intrinsically in the boundary integral
formulation, both along and perpendicular to the direction of
periodicity. This means that the imposition of artificial
(non-reflecting) boundaries, as would be required for finite element
and finite difference approaches, is not necessary.

A related class of potential problems can be found in the study of
water waves. Preston \emph{et al.} have derived and analyzed a
boundary integral formulation for the associated Dirichlet problem
on a half-space \cite{MP08} and presented an efficient and
super-algebraically convergent Nystr\"{o}m method for its numerical
solution in \cite{MP11}. These two papers are a subset of a wider
body of work on boundary integral methods for scattering from rough
surfaces (usually concerned with the solution of the Helmholtz
equation), and references can be found within. Boundary integral
methods for transmission problems can be found in \cite{CS85,KM88,MR08,GH11,SC13} and references therein. Work in this area is
generally concentrated on the solution of the Helmholtz equation in
both the interior and (unbounded) exterior of a two or
three-dimensional bounded domain, given some regularity properties
on the interface between the interior and exterior domains. Boundary
integral methods for transmission problems between various
configurations of bounded domains are considered in \cite{HH12i,HH12ii,CTG12}. The transmission problem for the Laplace equation is
also considered explicitly in \cite{CS85}.

The paper is structured as follows. In the next section we describe
the governing transmission problem for the Laplace equation and show
that it has at most one solution. In section \ref{BIM} we
reformulate this problem as a second-kind integral equation on the
interface, making use of the solution to the periodic half-space
boundary value problem, which plays the role of our boundary data.
The Dirichlet-to-Neumann operator for the interface is also
introduced leading to the solution of a first-kind integral equation
for the interface Neumann data. Section \ref{nystrom} then describes
the numerical solution of these two integral equations via the
Nystr\"{o}m method with appropriate quadrature rules. In particular
we draw on the theory presented in \cite{MP11,KS93,KA97,RK99} to
design a scheme that will exhibit super-algebraic convergence
properties. In section \ref{numerics} we demonstrate this
convergence rate via some numerical experiments with a set of
parameters that are appropriate for our application of modelling the
potential at the interface between a thin film of oil and the
surrounding air.

\section{The transmission problem for the Laplace equation}

Denote $H=\{(x,y)\in\mathbb{R}^2|y>0\}$ and let
$h:\mathbb{R}\rightarrow\mathbb{R}$ define an interface at $y=h(x)$
dividing $H$ into $\Omega_1=\{(x,y)\in H|y<h(x)\}$ and
$\Omega_2=\{(x,y)\in H|y>h(x)\}$. Here $\Omega_1$ represents a thin
film of oil with boundary $\Gamma_1=\Gamma_I\cup\Gamma_0$, where
$\Gamma_I=\{(x,y)\in\mathbb{R}^2|y=h(x)\}$ and
$\Gamma_0=\{(x,y)\in\mathbb{R}^2|y=0\}$. The domain $\Omega_2$
represents the surrounding air and has boundary $\Gamma_2=\Gamma_I$.
We will assume that $h$ is positive and $C^{k}$ for some $k\geq2$.
Let $BC^k(\bar{\Omega}_{\alpha})$ denote the space of bounded
functions in $C^{k}(\Omega_{\alpha})$ that can be continuously
extended into the closure $\bar{\Omega}_{\alpha}$ with $\alpha=1,2$.
Consider the following transmission problem for the electric
potential $\phi_{\alpha}\in BC^2({\bar\Omega}_{\alpha})$ for
$\alpha=1,2$:

\begin{align}\label{bvp1}
&\triangle\phi_{\alpha}=\:0
&\mathrm{in}\hspace{3mm}\Omega_\alpha,\\\label{intfc1}
&\phi_1=\:\phi_2 &\mathrm{on}\hspace{3mm}\Gamma_I,\\\label{intfc2}
&\epsilon_1 \frac{\partial\phi_1}{\partial
\nu}=\:\displaystyle\:\epsilon_2\frac{\partial\phi_2}{\partial\nu}
&\mathrm{on}\hspace{3mm}\Gamma_I.
\end{align}
In addition we prescribe the boundary conditions
\begin{align}\label{bvp2}
&\phi_1(x,y)=f(x) &\mathrm{on}\hspace{5mm}\Gamma_0,\\\label{bvp3}
&\frac{\partial\phi_2}{\partial y}(x,y)\rightarrow0
&\mathrm{as}\hspace{5mm}y\rightarrow\infty.
\end{align}
\begin{figure}[ht]
\begin{center}
\epsfig{file=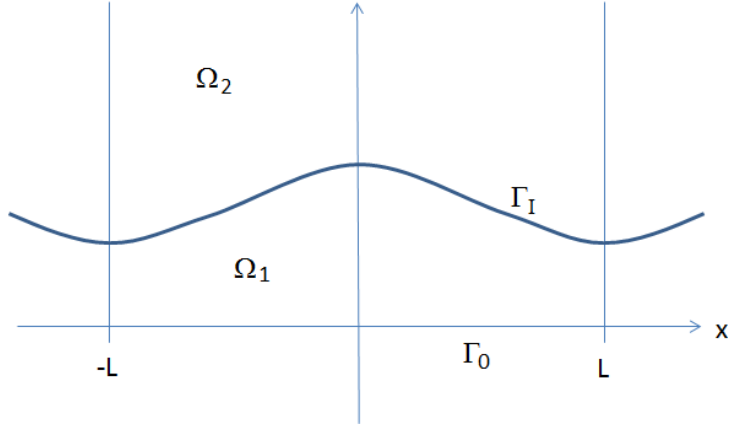, width=0.8\textwidth}
\end{center}
\caption{Problem setup}\label{prob1}
\end{figure}
Here $\epsilon_\alpha\in(0,\infty)$ is the dielectric constant in
$\Omega_\alpha$ and $\nu$ is the unit normal vector pointing out of
$\Omega_1$. Under the assumption that $h$ and $f$ are $2L-$periodic,
the problem is reduced to the study of a single periodic cell of $H$
with $-L\leq x<L$. Note that this is a reasonable assumption here
since $f$ represents the potential applied by a regularly spaced
array of
electrodes. The problem setup is shown in Figure \ref{prob1}.\\
\\
\begin{proposition}\label{P1}
Consider the transmission problem (\ref{bvp1}) to (\ref{bvp3}) with
$2L-$periodic boundary condition $f\in C(\mathbb{R})$ and
$2L-$periodic interface $0<h\in C^k(\mathbb{R})$ for some integer
$k\geq2$. If a solution $\phi_{\alpha}\in
BC^2({\bar\Omega}_{\alpha})$, $\alpha=1,2$ of this problem exists,
then it is unique.
\end{proposition}

\emph{Proof.} A similar argument to that given for (\cite{CS85},
Proposition 4.7) can be applied. We actually show that the boundary
value problem (\ref{bvp1}) to (\ref{bvp3}) with homogeneous boundary
condition $f=0$ only has the trivial solution $\phi_1=\phi_2=0$,
from which the proposition is a straightforward consequence.

Due to the periodicity restriction then a bound of the form $0<h<U$
holds for some constant $U$. Denoting the line $y=U$ by $\Gamma_U$
and applying Green's first identity on $\Omega_1$ and
$\Omega_2^*=\Omega_2\cap\{(x,y)\in\mathbb{R}^2|y<U\}$ gives
\begin{align}\label{O1}
\int_{\Gamma_0}\frac{\partial\phi_1}{\partial\nu}\phi_1dx+\int_{\Gamma_I}\frac{\partial\phi_1}{\partial\nu}\phi_1ds&=\int_{\Omega_1}|\nabla\phi_1|^2dA,\\\label{O2}
-\int_{\Gamma_U}\frac{\partial\phi_2}{\partial\nu}\phi_2dx-\int_{\Gamma_I}\frac{\partial\phi_2}{\partial\nu}\phi_2ds&=\int_{\Omega_2^*}|\nabla\phi_2|^2dA.
\end{align}

Note that periodicity means that the terms from integrating along
the vertical sides exactly cancel and the boundary condition $f=0$
means that the integral over $\Gamma_0$ vanishes. Combining
(\ref{O1}), (\ref{O2}) and the interface conditions (\ref{intfc1}),
(\ref{intfc2}) gives
\begin{align}\label{O3}
-\epsilon_2\int_{\Gamma_U}\frac{\partial\phi_2}{\partial\nu}\phi_2dx=\epsilon_2\int_{\Omega_2^*}|\nabla\phi_2|^2dA+\epsilon_1\int_{\Omega_1}|\nabla\phi_1|^2dA.
\end{align}
Taking the limit as $U\rightarrow\infty$ and using that
$\epsilon_\alpha>0$, $\alpha=1,2$ leads to the conclusion that
$\phi_\alpha$, $\alpha=1,2$ are both constant functions. The
boundary condition (\ref{bvp2}) and the continuity condition
(\ref{intfc1}) mean that these constants must both be zero. $\square$\\
\\

In the full multi-physics model, where the potential problem here is
coupled with a dynamic fluid-interface model, we wish to study the
quasi-time-dependent case where $f$ depends on time in the sense of
switching the applied potential on or off. This would mean solving
for $\phi_\alpha$ at a number of time-steps with varying interface
position $h$. In this work we simply consider modelling a static
potential for a range of different interface geometries $h$, with
the potential switched on. However we emphasize that the extension
of the methods developed here to the dynamic case would be
straightforward.

\section{Boundary integral formulation}\label{BIM}

In this section we recast the transmission problem
(\ref{bvp1})-(\ref{bvp3}) as a boundary integral equation on a
single periodic section of the interface $\Gamma_I$. From hereon we
make the abuse of notation that $\Gamma_I$, $\Gamma_0$ and
$\Omega_{\alpha}$ for $\alpha=1,\:2$ all refer to the restriction of
these curves / domains to a single periodic section from $x=-L$ to
$x=L$.

\subsection{Green's functions, layer potentials and the periodic half-plane solution}

A key ingredient in our boundary integral formulation will be the
$2L$-periodic half-plane Green's function for the Laplace equation
and its derivative with respect to $\nu$. The periodic Green's
function is given as (see for example \cite{CL99})
\begin{align}
G(\mathbf{x},\mathbf{x}_0)=-\frac{1}{2\pi}\ln\left(2\left|\sin\left(\frac{\pi}{2L}(z-z_0)\right)\right|\right),
\end{align}
where $\mathbf{x}=(x,y)\in H$, $z=x+iy$, $\mathbf{x}_0=(x_0,y_0)\in
\bar{H}$ and $z_0=x_0+iy_0$. It follows from the method of images that the
periodic Green's function on $H$ may be written
\begin{align}
G_H(\mathbf{x},\mathbf{x}_0)=G(\mathbf{x},\mathbf{x}_0)-G(\mathbf{x},\mathbf{x}_0'),
\end{align}
where $\mathbf{x}_0'=(x_0,-y_0)$ is the mirror image of
$\mathbf{x}_0$ in $\Gamma_0$. Computing the partial derivatives of
$G$ we obtain
\begin{align}\label{pd1}
\frac{\partial G_H}{\partial
x_0}(\mathbf{x},\mathbf{x}_0)&=\frac{1}{4L}\Re\left(\cot\left(\frac{\pi}{2L}(z-z_0)\right)-\cot\left(\frac{\pi}{2L}(z-z_0')\right)\right),\\\label{pd2}
\frac{\partial G_H}{\partial
y_0}(\mathbf{x},\mathbf{x}_0)&=-\frac{1}{4L}\Im\left(\cot\left(\frac{\pi}{2L}(z-z_0)\right)+\cot\left(\frac{\pi}{2L}(z-z_0')\right)\right),
\end{align}
where $z'_0=x_0-iy_0$. It is now straightforward from (\ref{pd1})
and (\ref{pd2}) to compute the derivative of $G_H$ with respect to
$\nu_0$, the unit normal vector at $\mathbf{x}_0$ pointing out of
$\Omega_1$.

It is convenient here to introduce interface single and double layer
potential operators $\mathcal{S}$ and $\mathcal{D}$, respectively,
acting on a continuous function $\phi$ at
$\mathbf{x}\in\Omega_\alpha$ ($\alpha=1,2$)
\begin{align}
\mathcal{S}\phi(\mathbf{x})&=\int_{\Gamma_I}G_H(\mathbf{x},\mathbf{x}_0)\phi(\mathbf{x}_0)d\Gamma(\mathbf{x}_0),\\
\mathcal{D}\phi(\mathbf{x})&=\int_{\Gamma_I}\frac{\partial
G_H}{\partial
\nu_0}(\mathbf{x},\mathbf{x}_0)\phi(\mathbf{x}_0)d\Gamma(\mathbf{x}_0).
\end{align}
In order to distinguish the case when $\mathbf{x}\in\Gamma_I$, we
rename the interface single and double layer potentials $V$ and $K$,
respectively.

Now let us consider the transmission problem
(\ref{bvp1})-(\ref{bvp3}) in the absence of an interface. The
problem is simplified to the study of a boundary value problem for
the Laplace equation on $\Omega_1=H$, with $2L$-periodic Dirichlet
data $f$ along $\Gamma_0$ (\ref{bvp2}), and condition (\ref{bvp3}).
Using Green's representation formula for $\mathbf{x}\in H$, the
solution of this half-plane problem $\phi_H$ may be expressed in the
form
\begin{align}\label{phiH1}
\phi_H(\mathbf{x})&=\int_{\Gamma_0}\left(G_H(\mathbf{x},\mathbf{x}_0)\frac{\partial
\phi_H}{\partial \nu_0}(\mathbf{x}_0)-\frac{\partial G_H}{\partial
\nu_0}(\mathbf{x},\mathbf{x}_0)f(x_0)\right)dx_0,\\\label{phiH2}
&=\int_{\Gamma_0}\frac{\partial G_H}{\partial
y_0}(\mathbf{x},\mathbf{x}_0)f(x_0)dx_0.
\end{align}
The uniqueness of $\phi_H$ follows from Prop. \ref{P1} with
$\epsilon_1=\epsilon_2$. We will make use of the solution $\phi_H$
in the next section when deriving the boundary integral formulation
for the transmission problem (\ref{bvp1})-(\ref{bvp3}).

\subsection{Boundary integral formulation for the transmission
problem}\label{BITrans}

We derive a direct boundary integral formulation for the
transmission problem (\ref{bvp1})-(\ref{bvp3}) using Green's
representation formula. This has the advantage, in comparison with
indirect formulations, that the unknowns appearing in the boundary
integral equations are the physical quantities (potentials) that we
wish to compute. We note also that a single layer potential solution
ansatz would fail here due to $G_H$ vanishing on $\Gamma_0$ where
the boundary condition $f$ is prescribed. A double layer potential
solution ansatz would lead to a more complicated equation for
computing the Neumann data on the interface than the method proposed
here. A further advantage of our direct approach is that it allows
our equations to be formulated only on the interface boundary
$\Gamma_I$, by making use of the half-space solution $\phi_H$,  as
described below. This advantage not only provides efficiency savings
though reducing the size of the domain to be discretized, but also
means that problems due to near-singularities associated with
boundary integral methods in long slender domains are avoided too.

We first apply Green's representation formula in $\Omega_2$ to give:
\begin{align}\label{domain2}
\phi_2=\mathcal{D}\phi_2-\mathcal{S}\left(\frac{\partial
\phi_2}{\partial\nu_0}\right).
\end{align}
Note that the integrals as $y\rightarrow\infty$ vanish due to both
$\phi_2$ and $G_H$ remaining bounded as $y\rightarrow\infty$ and
satisfying the radiation condition (\ref{bvp3}). Taking the limit as
the solution point tends to $\Gamma_I$ and applying the interface
conditions (\ref{intfc1}) and (\ref{intfc2}) yields
\begin{align}\label{O2eq}
V\left(\frac{\partial
\phi_1}{\partial\nu_0}\right)=\frac{\epsilon_2}{\epsilon_1}\left(-\frac{I}{2}+K\right)\phi_1.
\end{align}
The next step is to apply Green's representation formula in
$\Omega_1$ giving rise to
\begin{align}
\phi_1(\mathbf{x})=\int_{\Gamma_1}\left(G_H(\mathbf{x},\mathbf{x}_0)\frac{\partial
\phi_1}{\partial \nu_0}(\mathbf{x}_0)-\frac{\partial G_H}{\partial
\nu_0}(\mathbf{x},\mathbf{x}_0)\phi_1(\mathbf{x}_0)\right)d\Gamma_1(\mathbf{x}_0).
\end{align}
Splitting the integrals into the sum of an integral over $\Gamma_I$
and an integral over $\Gamma_0$ yields
\begin{align}\label{O1_2int}
\phi_1(\mathbf{x})=\mathcal{S}\left(\frac{\partial
\phi_1}{\partial\nu_0}\right)-\mathcal{D}\phi_1+\int_{\Gamma_0}\left(G_H(\mathbf{x},\mathbf{x}_0)\frac{\partial
\phi_1}{\partial \nu_0}(\mathbf{x}_0)-\frac{\partial G_H}{\partial
\nu_0}(\mathbf{x},\mathbf{x}_0)\phi_1(\mathbf{x}_0)\right)dx_0.
\end{align}
Applying (\ref{bvp2}) and (\ref{phiH1}), then (\ref{O1_2int})
simplifies to
\begin{align}\label{domain1}
\phi_1=\mathcal{S}\left(\frac{\partial
\phi_1}{\partial\nu_0}\right)-\mathcal{D}\phi_1+\phi_H.
\end{align}
Taking the limit as the solution point tends to $\Gamma_I$ as before
and rearranging yields
\begin{align}\label{O1eq}
V\left(\frac{\partial
\phi_1}{\partial\nu_0}\right)=\left(\frac{I}{2}+K\right)\phi_1-\phi_H.
\end{align}
Combining equations (\ref{O2eq}) and (\ref{O1eq}) results in the
following second-kind Fredholm integral equation for
$\phi=\phi_1=\phi_2$ on $\Gamma_I$
\begin{align}\label{Feq}
\left(I-2\frac{(\epsilon_2-\epsilon_1)}{(\epsilon_1+\epsilon_2)}K\right)\phi=\frac{2\epsilon_1}{(\epsilon_1+\epsilon_2)}\phi_H.
\end{align}
Clearly the case $\epsilon_1=\epsilon_2$ reduces to $\phi=\phi_H$ as
expected. Note that here we are effectively treating $\phi_H$ as our
boundary data on $\Gamma_I$, and that since $\phi_H$ is harmonic, it
is analytic. In the examples considered later, a simple closed form
expression will be available for $\phi_H$. However, in general it
may be necessary to approximate $\phi_H$ by, for example, a
truncated Fourier series (if solving the half-plane problem via
separation of variables) or quadrature (if computing $\phi_H$
directly from the boundary integral formula (\ref{phiH2})).

Let us now consider the operator $K$ in further detail. Since for
$\mathbf{x}\in\Gamma_I$ then $(x,y)=(x,h(x))$, then if we also have $\mathbf{x}_0\in\Gamma_I$ we may write
\begin{align}
T_H(x,x_0):=\frac{\partial
G_H}{\partial\nu_0}(\mathbf{x},\mathbf{x}_0),
\end{align}
and hence, with a slight abuse of notation,
\begin{align}
K\phi(x)=\int_{-L}^L T_H(x,x_0)\phi(x_0)\sqrt{1+h'(x_0)^2}d x_0.
\end{align}
Considering therefore a single periodic strip with $x\in[-L,L]$,
then the behavior of $K$ close to $x=x_0$ is identical to that for
the double layer potential on a simple closed curve studied in
(\cite{KA97}, Chap. 7). This is evident from taking limits in
(\ref{pd1}) and (\ref{pd2}) as $(x-x_0)\rightarrow0$. Since $h\in
C^k([-L,L])$ for some $k\geq2$, then it follows from (\cite{KA97},
Chap. 7) that
\begin{align}\label{ksmooth}
T_H\in C^{k-2}([-L,L]\times[-L,L]).
\end{align}
As a result of this $K$ is a compact operator on $C([-L,L])$ with
the maximum norm (see for example, \cite{RK99}, Theorem 2.21). Using
the theory of compact integral operators \cite{RK99} we now prove
the existence and uniqueness of solutions to the integral equation
(\ref{Feq}). We first prove the following
\begin{align}\label{DoubImp}
\begin{split}
f=0\: & \Leftrightarrow\:\phi_H=0\\
& \Leftrightarrow\:\phi=0. \end{split}
\end{align}
The implications $f=0\Rightarrow\phi_H=0$ and $f=0\Rightarrow\phi=0$
are due to the uniqueness of solution to both the transmission
problem (Prop. \ref{P1}) and the periodic half-plane problem, which
follows from Prop. \ref{P1} with $\epsilon_1=\epsilon_2$. The
implication $\phi_H=0\Leftarrow\phi=0$ follows from (\ref{Feq}).
Finally, $f=0\Leftarrow\phi_H=0$ relies on the fact that $\phi_H$
can be continuously extended to $\Gamma_0$ (\cite{MP08}, Theorem
2.3) so that
$$\lim_{y\rightarrow0}\phi_H(\mathbf{x})=f(x),$$
which may be derived from the jump properties of the double layer
potential (\ref{phiH2}).

Hence assuming $\phi_H=0$, then (\ref{DoubImp}) gives that $\phi=0$
is the only solution of
\begin{align}\label{Feq0}
\left(I-2\frac{(\epsilon_2-\epsilon_1)}{(\epsilon_1+\epsilon_2)}K\right)\phi=0.
\end{align}
Hence the operator
\begin{align}
\left(I-2\frac{(\epsilon_2-\epsilon_1)}{(\epsilon_1+\epsilon_2)}K\right)
\end{align}
is injective on $C([-L,L])$. It now follows from (\cite{RK99}, Thm.
3.4) that (\ref{Feq}) is uniquely solvable and that the solution
depends continuously on $\phi_H$.
\begin{theorem}
The integral operator $(I-2\mu K):C([-L,L])\rightarrow C([-L,L])$ is
invertible with bounded inverse for any $|\mu|<1$. Consequently,
there exists a unique solution $\phi_{\alpha}\in
BC^2({\bar\Omega}_{\alpha})$ to the transmission problem
(\ref{bvp1}) to (\ref{bvp3}) considered in Proposition 1.
\end{theorem}

\subsection{The Dirichlet-to-Neumann (DtN) operator}

For our application in tracking the dynamic evolution of a thin film
of oil, the boundary data for the fluid equations will depend on the
normal and tangential derivatives of $\phi$ on $\Gamma_I$.
Furthermore the normal derivative would be needed for computing the
domain potentials via (\ref{domain2}) and (\ref{domain1}). Later we
discuss interpolation formulae for obtaining the tangential
derivative. The normal derivative will be computed using the DtN
operator. Combining (\ref{intfc2}) and (\ref{O2eq}) yields the
following first-kind integral equation for $\partial
\phi_\alpha/\partial\nu_0$, $\alpha=1,2$:
\begin{align}\label{dtn1}
V\left(\frac{\partial
\phi_\alpha}{\partial\nu_0}\right)=\frac{\epsilon_{2}}{\epsilon_\alpha}\left(-\frac{I}{2}+K\right)\phi.
\end{align}
The DtN operator will be bounded as a map from $H^{1/2}(\Gamma_I)$
to $H^{-1/2}(\Gamma_I)$ if we can prove the invertibility of
$V:H^{-1/2}(\Gamma_I)\rightarrow H^{1/2}(\Gamma_I)$ with bounded
inverse, where $H^{\pm1/2}$ are the usual Sobolev spaces
(\cite{RK99}, Sec. 8.2).  To prove this we need to show that $V$ is
bounded and elliptic so the Lax-Milgram theorem gives that
$V^{-1}:H^{1/2}(\Gamma_I)\rightarrow H^{-1/2}(\Gamma_I)$ is bounded.
The proof that $V:H^{-1/2}(\Gamma_I)\rightarrow H^{1/2}(\Gamma_I)$
is bounded can be done using a similar argument to the one presented
in (\cite{RK99}, Theorem 8.23) for the single layer potential on a
simple closed $C^k$ curve in $\mathbb{R}^2$. The proof of
ellipticity is again similar to that for the case of a closed curve
in $\mathbb{R}^2$. A rough outline restricted to
$V:H_{*}^{-1/2}(\Gamma_I)\rightarrow H^{1/2}(\Gamma_I)$, where
$$H^{-1/2}_{*}(\Gamma_I)= \left\{\tilde{\sigma}\in
H^{-1/2}(\Gamma_I)\:\left|\int_{\Gamma_{I}}\tilde{\sigma}ds=0\right.\right\},$$
is presented below.

Let $\phi=\mathcal{S}\sigma$, where
$\phi_\alpha=\phi|_{\Omega_\alpha}$ for $\alpha=1,2$ and $\sigma$ is
a so-called layer density to be determined. The definition of
$\mathcal{S}$ allows us to deduce that $\phi|_{\Gamma_0}=0$ and
\begin{equation}\label{PhiLim}
\lim_{y\rightarrow\infty}\frac{\partial\phi_2}{\partial y}=
0.
\end{equation}
Furthermore, $V$ and its normal derivative will satisfy the standard
jump relations associated with the classical single layer potential
operator for the Laplace equation in two dimensions. This is simply
a consequence of the fact that $G_H$ is constructed using the method
of images to both incorporate periodicity and to restrict to a
half-plane. The resulting additional contributions to the classical
operator and their normal derivatives will all be continuous across
$\Gamma_I$ (on the periodic section considered). Hence on $\Gamma_I$
we have that $\phi_1=\phi_2$ and
\begin{equation}\label{Sjump}
\frac{\partial\phi_1}{\partial \nu}-\frac{\partial\phi_2}{\partial
\nu}=\sigma.
\end{equation}
From the above it is then straightforward to deduce that
\begin{equation}\label{Sjumpvar}
\int_{\Gamma_I}\left(\phi_1\frac{\partial\phi_1}{\partial
\nu}-\phi_2\frac{\partial\phi_2}{\partial \nu}\right)
ds=\int_{\Gamma_I}\sigma V\sigma ds.
\end{equation}
Arguing as in the proof of Prop. \ref{P1} using Green's first
identity, (\ref{PhiLim}) and that $\phi|_{\Gamma_0}=0$ leads to
\begin{align}\label{green1}
\int_{\Gamma_I}\phi_\alpha\frac{\partial\phi_\alpha}{\partial\nu}ds=(-1)^{\alpha-1}\int_{\Omega_\alpha}|\nabla\phi_{\alpha}|^2dA.
\end{align}
Combining (\ref{Sjumpvar}) and (\ref{green1}) we find
\begin{align}\label{p1}
\begin{split}
\int_{\Gamma_I}\sigma V\sigma
ds & =\int_{\Omega_1}|\nabla\phi_{1}|^2dA+\int_{\Omega_2}|\nabla\phi_{2}|^2dA\vspace{2mm}\\
&=|\phi_{1}|^2_{H^1(\Omega_1)}+|\phi_{2}|^2_{H^1(\Omega_2)}.
\end{split}
\end{align}
Now, considering the norm of the layer density we have
\begin{align}\label{p2}
\begin{split}
\|\sigma\|^2_{H^{-1/2}(\Gamma_I)} & =\left\|\frac{\partial\phi_1}{\partial \nu}-\frac{\partial\phi_2}{\partial
\nu}\right\|^2_{H^{-1/2}(\Gamma_I)}\\
&\leq 2\left(\left\|\frac{\partial\phi_1}{\partial \nu}\right\|^2_{H^{-1/2}(\Gamma_I)}+\left\|\frac{\partial\phi_2}{\partial
\nu}\right\|^2_{H^{-1/2}(\Gamma_I)}\right)\\
& \leq C\left(|\phi_{1}|^2_{H^1(\Omega_1)}+|\phi_{2}|^2_{H^1(\Omega_2)}\right).
\end{split}
\end{align}
The final inequality follows by applying the Cauchy-Schwarz
inequality in the first Green identity  to get
\begin{align}
\left|\int_{\Gamma_I}\psi\frac{\partial\phi_\alpha}{\partial\nu}ds\right|\leq
\:|\phi_{\alpha}|_{H^1(\Omega_\alpha)}|\psi|_{H^{1}(\Omega_\alpha)},
\end{align}
with $\psi=\mathcal{S}\tilde{\sigma}$ for some $\tilde{\sigma}\in
H_{*}^{-1/2}(\Gamma_{I})$. Then applying the inverse trace theorem
to $\psi|_{\Gamma_I}\in H^{1/2}(\Gamma_I)$ and using duality of
$H^{\pm1/2}(\Gamma_I)$ leads to an inequality of the form
\begin{align}
\left\|\frac{\partial\phi_\alpha}{\partial\nu}\right\|_{H^{-1/2}(\Gamma_I)}=\sup_{\psi\in
H^{1/2}(\Gamma_I)\setminus\{0\}}\frac{\left|\int_{\Gamma_I}\psi\frac{\partial\phi_\alpha}{\partial\nu}ds\right|}{\|\psi\|_{H^{1/2}(\Gamma_I)}}\leq
C_{\alpha}\:|\phi_{\alpha}|_{H^1(\Omega_\alpha)},
\end{align}
where $C\:\mathrm{and}\:C_{\alpha}$ are positive constants with
$C=2\max (C_\alpha)^2$. To conclude the argument we combine
(\ref{p1}) and (\ref{p2}) to give ellipticity, that is
\begin{align}
C\int_{\Gamma_I}\sigma V\sigma
ds\geq\|\sigma\|^2_{H^{-1/2}(\Gamma_I)}.
\end{align}

\section{Discretization using the Nystr\"{o}m method}\label{nystrom}

In this section we describe the discretization of the integral
equations (\ref{Feq}) and (\ref{dtn1}) using the Nystr\"{o}m method
with suitable quadratures. In particular, we draw on the theory
presented in \cite{MP11,KS93,KA97,RK99} in order to design a method
with super-algebraic convergence, as demonstrated by the numerical
results in the next section. It should be noted that to obtain these
convergence rates requires that $h\in C^{\infty}([-L,L])$, and so
from hereon we make this assumption. In the application to modelling
the potential in a layer of oil, where the (time-dependent) oil-air
interface position will only be given as a set of equi-spaced
coordinates, we will use interpolation by trigonometric polynomials
to ensure the required smoothness. This is a good choice of
interpolation scheme due to the periodicity of $h$, and the
availability of a highly efficient FFT based implementation. Even
though $h$ will be given explicitly in the examples here, we will
interpolate anyway to ensure the algorithm can be applied more
generally.

\subsection{Discretization for a smooth kernel}

Let us first consider the second-kind equation (\ref{Feq}). Under
the assumption of an infinitely differentiable interface, we also
have an infinitely differentiable kernel and as discussed before,
the data term $\phi_H$ is also infinitely differentiable. In this
situation a simple application of the trapezoidal rule yields a
super-algebraically convergent method \cite{MP11,KA97}. To
implement this scheme we first use that $\mathbf{x}=(x,h(x))$ on
$\Gamma_I$, and write integrals over $\Gamma_I$ in the form
\begin{align}
\int_{\Gamma_I}F(\textbf{x})d\Gamma(\textbf{x})=\int_{-L}^{L}\tilde{F}(x)\sqrt{1+h'(x)^2}dx.
\end{align}
Applying the trapezoidal rule with $n$ subintervals gives
\begin{align}\label{Trap}
\int_{\Gamma_I}F(\textbf{x})d\Gamma(\textbf{x})\approx\frac{2L}{n}\sum_{j=1}^{n}\tilde{F}\left(x_j\right)\sqrt{1+h'\left(x_j\right)^2},
\end{align}
where $x_j=-L+2L(j-1)/n$. Note that because we have assumed above
that $h$ is given by trigonometric interpolation, its derivative may
be computed simply by differentiating its Fourier components as
described in \cite{MP11}.

Applying the formula (\ref{Trap}) to the integral in the definition
of $K$ yields the following approximation
\begin{align}\label{TrapK}
K\phi(x_i,\:h(x_i))\approx\frac{2L}{n}\sum_{j=1}^{n}k_{i,j}\phi(x_j)\sqrt{1+h'\left(x_j\right)^2},
\end{align}
where
\begin{align}\label{Kij1}
k_{i,j}=T_H(x_i,x_j).
\end{align}
The following explicit formula for $k_{j,j}$ may be derived (for
similar derivations see \cite{MP11} and \cite{KA97}, Ch. 7):
\begin{align}\label{Kij2}
k_{j,j}=\frac{h''(x_j)}{4\pi(1+h'(x_j)^2)}+\frac{1}{4L\sqrt{1+h'\left(x_j\right)^2}}\coth\left(\frac{\pi
}{L}h(x_j)\right).
\end{align}
As with the first derivative, $h''$ may be computed simply by
differentiating the Fourier components for $h$ twice. For $i\neq j$
we can use the formula obtained directly from (\ref{pd1}) and
(\ref{pd2}), that is,
\begin{align*}
k_{i,j}=\frac{1}{4L\sqrt{1+h'(x_j)^2}}&\left(h'(x_j)\:\Re\left(\cot\left(\frac{\pi}{2L}(z_i-z_j')\right)-\cot\left(\frac{\pi}{2L}(z_i-z_j)\right)\right)\right.\vspace{2mm}\\
&\left.-\Im\left(\cot\left(\frac{\pi}{2L}(z_i-z_j)\right)+\cot\left(\frac{\pi}{2L}(z_i-z_j')\right)\right)\right).
\end{align*}
Here we have denoted $z_j=x_j+ih(x_j)$ and $z_j'=x_j-ih(x_j)$.
Notice that (\ref{Kij2}) illustrates the smoothness result
(\ref{ksmooth}). Applying the approximation (\ref{TrapK}) to the
second-kind integral equation (\ref{Feq}) leads to the following
approximate Nystr\"{o}m scheme for the approximate solution $\phi^n$
\begin{align}\label{Nyst2ndKind}
\phi^{n}(x_i)+\frac{4L(\epsilon_1-\epsilon_2)}{n(\epsilon_1+\epsilon_2)}\sum_{j=1}^{n}k_{i,j}\phi^{n}(x_j)\sqrt{1+h'\left(x_j\right)^2}=\frac{2\epsilon_1}{(\epsilon_1+\epsilon_2)}\phi_H(x_i,h(x_i))
\end{align}
for $i=1,...,n$. The super-algebraic convergence of $\phi^n$ to
$\phi$ for increasing $n$ is a consequence of (\cite{MP11}, Thm.
3.12).

\subsection{Discretization for a kernel with a logarithmic singularity}

In this section we consider a Nystr\"{o}m method for the solution of
equation (\ref{dtn1}), which is a first-kind integral equation for
$\partial \phi_\alpha/\partial\nu_0$, $\alpha=1,2$. In particular we
note that the kernel function $G_H$ of the operator $V$ contains a
logarithmic singularity and may be written in the form
\begin{align}\label{GHKressSloan1}
G_H(\mathbf{x},\mathbf{x}_0)=-\frac{1}{4\pi}\ln\left(4\sin^2\left(\frac{\pi}{2L}(x-x_0)\right)\right)+\tilde{G}(x,x_0)-G(\mathbf{x},\mathbf{x}_0'),
\end{align}
with
\begin{align*}
\tilde{G}(x,x_0)=-\frac{1}{4\pi}\left\{\ln\left(4\left(\sin^2\left(\frac{\pi}{2L}(x-x_0)\right)+\sinh^2\left(\frac{\pi}{2L}(h(x)-h(x_0))\right)\right)\right)\right.\\
\left.-\ln\left(4\sin^2\left(\frac{\pi}{2L}(x-x_0)\right)\right)\right\}.
\end{align*}
On $\Gamma_I$ with $y=h(x)$ for $h\in C^{\infty}([-L,L])$, then both
$\tilde{G}(x,x_0)$ and $G(\mathbf{x},\mathbf{x}_0')$ are infinitely
differentiable with respect to $x, x_0\in[-L,L]$. In particular, the
diagonal term for $\tilde{G}$ is given by
\begin{align}\label{GHKressSloan3}
\tilde{G}(x,x)=-\frac{1}{4\pi}\ln\left(1+h'(x)^2\right).
\end{align}
Note that (\ref{GHKressSloan1}) implies that the first-kind equation
(\ref{dtn1}) falls into the class of first kind equations analyzed
in \cite{KS93}. We therefore base our quadrature method here on the
approach suggested in \cite{KS93} (see also \cite{RK99} and
references therein). Super-algebraic convergence rates will then be
attained due to (\cite{KS93}, Thm. 2.3) if the right hand side of
(\ref{dtn1}) is infinitely differentiable. Note however that since
in general only an approximate numerical solution of (\ref{Feq})
will be available, then $\phi$ will need to be replaced by $\phi^n$
in (\ref{dtn1}). In order to satisfy the necessary criteria for the
super-algebraic convergence rates predicted by (\cite{KS93}, Thm.
2.3), we may interpolate $\phi^n$ using trigonometric polynomials to
obtain an infinitely differentiable interpolant. The right hand side
of (\ref{dtn1}) may then be written
$\epsilon_{2}(-\phi^n+2K\phi^n)/(2\epsilon_\alpha)$, which is
clearly infinitely differentiable since
\begin{align}
\frac{d^m}{dx^m}(K\phi^n)(\mathbf{x})=\int_{-L}^L\frac{\partial^m}{\partial
x^m}T_H(x,x_0)\phi^n(x_0)\sqrt{1+h'(x_0)^2}dx_0,
\end{align}
and $T_H$  is infinitely differentiable for $h\in
C^{\infty}([-L,L])$.

We now outline the quadrature rule we employ to approximate the
interface single layer potential $V$. For the term containing a
logarithmic singularity we employ a quadrature rule of the form
\begin{align}\label{FQuad1}
-\frac{1}{4\pi}\int_{-L}^{L}\ln\left(4\sin^2\left(\frac{\pi}{2L}(x-x_0)\right)\right)F(x_0)dx_0\approx\sum_{j=1}^{2N}R_{j}(x)F\left(x_j\right),
\end{align}
for positive integer $N=n/2$ (assuming $n$ is even), with
$x_j=-L+L(j-1)/n$, $j=1,...,n$, as before. The quadrature weight
function $R_{j}(x)$ is given by
\begin{align}\label{FQuad2}
R_{j}(x)=\frac{L}{2\pi
n}\left\{\sum_{m=1}^{N-1}\frac{1}{m}\cos\left(\frac{m\pi}{L}(x-x_j)\right)+\frac{1}{n}\cos\left(\frac{N\pi}{L}(x-x_j)\right)\right\}.
\end{align}
This choice of quadrature computes the integral in (\ref{FQuad1})
exactly when $F$ has been replaced by its trigonometric
interpolation polynomial. To see this replace $F$ in (\ref{FQuad1})
by the Lagrange trigonometric polynomial of order $j$, then the
formula (\ref{FQuad2}) may be derived for the integral on the left
hand side; see \cite{KS93} and (\cite{RK99}, pp. 208) for details.

The remaining integrals involved in approximating $V$ have smooth
kernels and thus may be well-approximated using the trapezoidal rule
as before. We therefore arrive at the following approximation for
the interface single layer potential
\begin{align}\label{slpapprox}
V\frac{\partial\phi_\alpha}{\partial\nu_0}(x_i,\:h(x_i))\approx\sum_{j=1}^{n}\left(R_{j}(x_i)+\frac{2L}{n}s_{i,j}\right)\frac{\partial\phi_\alpha}{\partial\nu_0}(x_j)\sqrt{1+h'(x_j)^2}
\end{align}
for $\alpha=1,2$ and $i=1,...,n$. Here $s_{i,j}$ is the smooth part
of the integrand for $V$ given by
\begin{align*}
s_{i,j}=-\frac{1}{4\pi}\left\{\ln\left(4\left(\sin^2\left(\frac{\pi}{2L}(x_i-x_j)\right)+\sinh^2\left(\frac{\pi}{2L}(h(x_i)-h(x_j))\right)\right)\right)\right.\\
-\ln\left(4\left(\sin^2\left(\frac{\pi}{2L}(x_i-x_j)\right)+\sinh^2\left(\frac{\pi}{2L}(h(x_i)+h(x_j))\right)\right)\right)\\
\left.-\ln\left(4\sin^2\left(\frac{\pi}{2L}(x_i-x_j)\right)\right)\right\},
\end{align*}
for $i\neq j$ and
\begin{align}
s_{j,j}=-\frac{1}{4\pi}\left\{\ln\left(1+h'(x_j)^2\right)-\ln\left(4\sinh^2\left(\frac{\pi}{L}h(x_j)\right)\right)\right\}.
\end{align}
Applying the approximation (\ref{slpapprox}) to the first-kind
integral equation (\ref{dtn1}) leads to the following Nystr\"{o}m
scheme for the approximate solution $\partial
\phi^n_\alpha/\partial\nu_0$
\begin{align}\label{Nyst1stKind}
\begin{array}{l}
\displaystyle
\sum_{j=1}^{n}\left(R_{j}(x_i)+\frac{2L}{n}s_{i,j}\right)\frac{\partial\phi^n_\alpha}{\partial\nu_0}(x_j)\sqrt{1+h'(x_j)^2}\\
\displaystyle\hspace{25mm}=\frac{\epsilon_2}{\epsilon_\alpha}\left(\frac{-\phi^n(x_i)}{2}+\frac{2L}{n}\sum_{j=1}^{n}k_{i,j}\phi^{n}(x_j)\sqrt{1+h'\left(x_j\right)^2}\right),
\end{array}
\end{align}
for $i=1,...,n$ and $\alpha=1,2$. Here $\phi^n$ and $k_{i,j}$ are
computed as described in the previous section.

\section{Numerical experiments}\label{numerics}
We test our boundary integral model and numerical solution scheme
using the same boundary condition for (\ref{bvp2}) as employed in
\cite{CB11} for a static applied potential $f(x)=A\cos(\pi x/L)$,
where $A$ is a constant amplitude term. Note that one can then
derive the half-space solution $\phi_H$, either using separation of
variables or directly from the boundary integral formula
(\ref{phiH2}), to get
\begin{align}\label{phiH}
\phi_{H}(\mathbf{x})=A\cos(\pi x/L)\exp(-\pi y/L).
\end{align}
In addition we set the parameter $\epsilon_1=8$ and $\epsilon_2=1$
to reflect the values for oil and air in our application of interest
\cite{CB11}. For the interface position $h(x)$ we consider three
possibilities:
\begin{enumerate}[{I}1]
\item A constant profile $h(x)=h_0$ since in this case the results can
be checked against an exact solution.
\item A single period sine curve
$h(x)=h_0(1+0.2\sin(\pi x/L))$ of mean height $h_0$, similar to that
used for numerical tests in \cite{MP11}.
\item A double period cosine curve $h(x)=h_0(1-\cos(2\pi x/L)/(2\pi))$ as considered in
Sec. II of \cite{CB11}, which is representative of a typical
interface geometry that would be encountered in application to an
electrified oil-air interface.
\end{enumerate}
Note that to allow the input of numeric interface position data,
trigonometric interpolation is used to represent the interface
position (exactly for the examples above), and derivatives of $h$
are computed via differentiation of its Fourier components as
described before.

We use the methods described in the previous section to approximate
the interface potential $\phi$ by $\phi^n$ and the normal derivative
$\partial\phi_{\alpha}/\partial\nu$ by
$\partial\phi_{\alpha}^n/\partial\nu$ ($\alpha=1,2$). We also
approximate the tangential derivative
$\partial\phi_{\alpha}/\partial\tau$, where $\tau$ denotes the
tangent vector to $\Gamma_I$. To do this we first compute the
derivative with respect to $x$ by applying trigonometric
interpolation to $\phi^n$ using an FFT and then differentiating the
Fourier components (exactly the same as the procedure we carry out
for the interface position function $h$). One then obtains the
tangential derivative by correcting for arc-length via division by a
factor of $\sqrt{1+h'(x)^2}$. A similar procedure is described and
rigorously analyzed in \cite{MP11}, where it is shown
super-algebraic convergence will also be achieved in the
approximation of $\partial\phi_{\alpha}/\partial\tau$ by
$\partial\phi^{n}_{\alpha}/\partial\tau$, $\alpha=1,2$. Note that
here the continuity of $\partial\phi^{n}_{\alpha}/\partial\tau$
across $\Gamma_I$ follows from the continuity of $\phi^{n}$ across
$\Gamma_I$, so we may simply write $\partial\phi^{n}/\partial\tau$
for both $\partial\phi^{n}_{\alpha}/\partial\tau$, $\alpha=1,2$.

\begin{figure}[ht]
\begin{center}
\epsfig{file=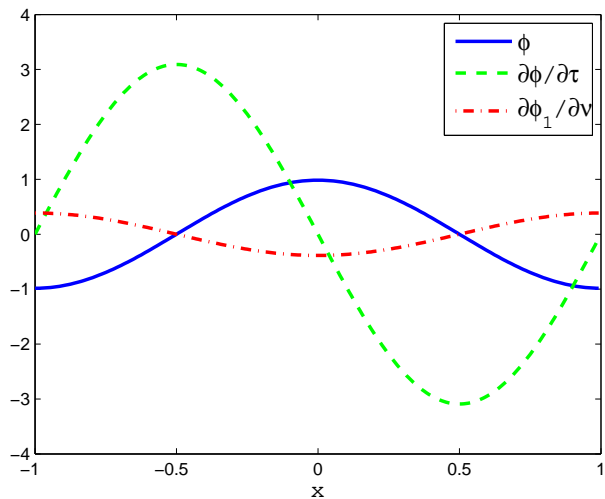, width=0.8\textwidth}
\end{center}
\caption{Analytic solution of the potential and its tangential and
normal derivatives on a flat interface.}\label{ExactFlat}
\end{figure}

In the case of interface I1, an analytic solution to the
transmission problem (\ref{bvp1}) to (\ref{bvp3}) with boundary data
$f(x)=A\cos(\pi x/L)$ can be derived using separation of variables
to give
\begin{align*}
\phi_1(\mathbf{x})&=\frac{A\cos(\pi
x/L)\left((\epsilon_1-\epsilon_2)\exp(\pi
y/L)+(\epsilon_1+\epsilon_2)\exp((2h_0-y)\pi/L)\right)}{(\epsilon_1+\epsilon_2)\exp(2\pi
h_0/L)+(\epsilon_1-\epsilon_2)},\\
\phi_2(\mathbf{x})&=\frac{2A\epsilon_1\cos(\pi
x/L)\exp((2h_0-y)\pi/L)}{(\epsilon_1+\epsilon_2)\exp(2\pi
h_0/L)+(\epsilon_1-\epsilon_2)}.
\end{align*}
From the solution above it is easy to calculate $\phi$ and its
tangential and normal derivatives, which here correspond with
partial derivatives with respect to $x$ and $y$, for interface I1. A
plot of $\phi$ for $A=L=1$ and a thin film with $h_0=0.03L$
(comparable to the parameter choice of \cite{CB11}) is given in Fig.
\ref{ExactFlat}. The tangential and normal derivatives are also
shown.

\begin{figure}[ht]
\begin{center}
\epsfig{file=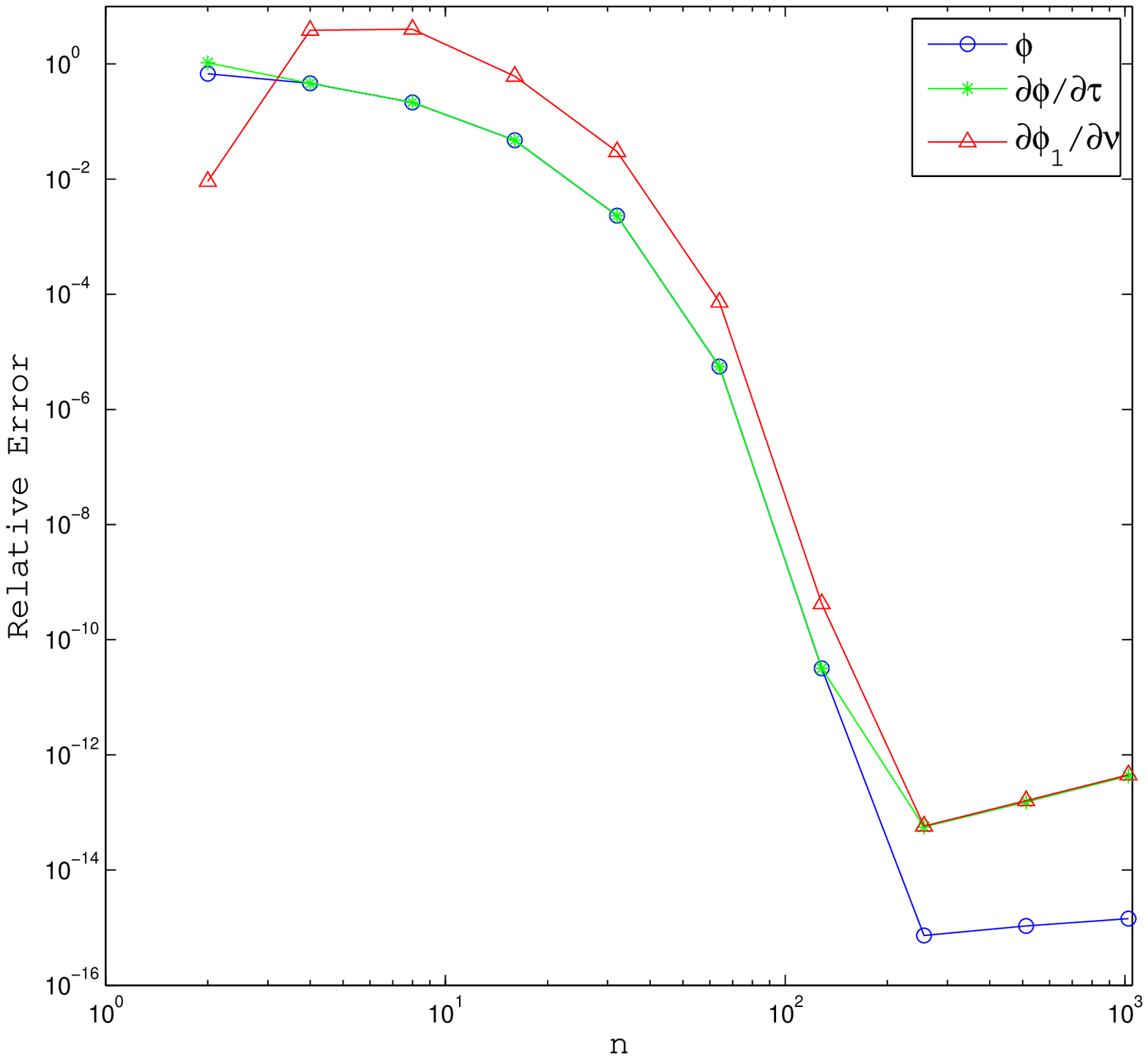, width=0.8\textwidth}
\end{center}
\caption{Relative errors in the potential and its tangential and
normal derivatives on a flat interface.}\label{FlatError}
\end{figure}

We use these analytic solutions to compute the discrete $l_2$ error
in our numerical approximations of $\phi$,
$\partial\phi/\partial\tau$ and $\partial\phi/\partial\nu$ with the
above choice of parameters. The results are shown in Fig.
\ref{FlatError} and demonstrate the predicted super-algebraic
convergence since the rate of convergence is speeding up as $n$ is
increased, until rounding errors become significant for $n>128$. The
results for $\phi$ and $\partial\phi/\partial\tau$ are almost
identical until the rounding errors dominate, and then the errors
for $\partial\phi/\partial\tau$ and $\partial\phi/\partial\nu$ are
very similar. The general trends are very similar to those shown the
Dirichlet problem studied in \cite{MP11}.

\begin{table}
\caption{Relative errors and estimated orders of convergence for
interface I1}
\begin{center} \footnotesize
\begin{tabular}{ccccccc} \hline
& \multicolumn{2}{c}{Results for $\phi^n$} & \multicolumn{2}{c}{Results for $\partial\phi^n/\partial\tau$} & \multicolumn{2}{c}{Results for $\partial\phi_1^n/\partial\nu$}\\
\hline
$n$  &  Error & EOC &  Error & EOC & Error & EOC\\
\hline
16  & \lower.3ex\hbox{4.506e-02}  & \lower.3ex\hbox{-}  & \lower.3ex\hbox{4.506e-02} & \lower.3ex\hbox{-}  & \lower.3ex\hbox{0.6025}     & \lower.3ex\hbox{-}  \\
32  & \lower.3ex\hbox{2.307e-03}  & \lower.3ex\hbox{4.29}  & \lower.3ex\hbox{2.307e-03} & \lower.3ex\hbox{4.29}  & \lower.3ex\hbox{2.983e-02}   & \lower.3ex\hbox{4.34}  \\
64  & \lower.3ex\hbox{5.551e-06}  & \lower.3ex\hbox{8.70}  & \lower.3ex\hbox{5.551e-06} & \lower.3ex\hbox{8.70}  & \lower.3ex\hbox{7.279e-05}   & \lower.3ex\hbox{8.68}  \\
128 & \lower.3ex\hbox{3.199e-11}  & \lower.3ex\hbox{17.40} & \lower.3ex\hbox{3.199e-11} & \lower.3ex\hbox{17.40} & \lower.3ex\hbox{4.227e-10}  & \lower.3ex\hbox{17.39} \\
256 & \lower.3ex\hbox{1.187e-15}  & \lower.3ex\hbox{14.72} & \lower.3ex\hbox{1.611e-13} & \lower.3ex\hbox{7.63}  & \lower.3ex\hbox{1.566e-13}  & \lower.3ex\hbox{11.40} \\
512 & \lower.3ex\hbox{1.769e-15}  & \lower.3ex\hbox{-0.58} & \lower.3ex\hbox{4.579e-13} & \lower.3ex\hbox{-1.51} & \lower.3ex\hbox{4.724e-13}  & \lower.3ex\hbox{-1.59} \\
\hline
\end{tabular}
\end{center}
\label{FlatTab}
\end{table}

Since analytic solutions are not available for the cases of
interface I2 and I3, we compute the discrete $l_2$ relative errors
and estimate their rates of convergence by using $\phi^{2n}$,
$\partial\phi^{2n}/\partial\tau$ and $\partial\phi^{2n}/\partial\nu$
in the role of the `exact' solution for $\phi^{n}$,
$\partial\phi^{n}/\partial\tau$ and $\partial\phi^{n}/\partial\nu$,
respectively (using only the coincident data points). For
consistency we also show the results for the case of interface I1 in
this form, which are presented in Table \ref{FlatTab}. The results
clearly reflect those shown in Fig. \ref{FlatError}. The estimated
order of convergence (EOC) is computed via
$\log_2(\mathrm{Error}(n)/\mathrm{Error}(n/2))$, where we are
referring to the discrete $l_2$ relative error as a function of $n$.
The super-algebraic convergence is clear from the increasing EOC,
which eventually reaches a value of around $17$ in all cases before
rounding errors begin to destroy the convergence rate when $n=256$.

\begin{table}
\caption{Relative errors and estimated orders of convergence for
interface I2}
\begin{center} \footnotesize
\begin{tabular}{ccccccc} \hline
& \multicolumn{2}{c}{Results for $\phi^n$} & \multicolumn{2}{c}{Results for $\partial\phi^n/\partial\tau$} & \multicolumn{2}{c}{Results for $\partial\phi_1^n/\partial\nu$}\\
\hline
$n$  &  Error & EOC &  Error & EOC & Error & EOC\\
\hline
16  & \lower.3ex\hbox{4.872e-02}  & \lower.3ex\hbox{-}  & \lower.3ex\hbox{5.462e-02} & \lower.3ex\hbox{-}  & \lower.3ex\hbox{7.231e-01}   & \lower.3ex\hbox{-}  \\
32  & \lower.3ex\hbox{3.248e-03}  & \lower.3ex\hbox{3.91}  & \lower.3ex\hbox{4.610e-03} & \lower.3ex\hbox{3.57}  & \lower.3ex\hbox{4.832e-02}   & \lower.3ex\hbox{3.90}  \\
64 & \lower.3ex\hbox{1.721e-05}  & \lower.3ex\hbox{7.56} & \lower.3ex\hbox{3.410e-05} & \lower.3ex\hbox{7.08} & \lower.3ex\hbox{2.806e-04}  & \lower.3ex\hbox{7.43} \\
128 & \lower.3ex\hbox{7.230e-10}  & \lower.3ex\hbox{14.54}  & \lower.3ex\hbox{2.055e-09} & \lower.3ex\hbox{14.02}  & \lower.3ex\hbox{1.293e-08}  & \lower.3ex\hbox{14.41} \\
256 & \lower.3ex\hbox{1.861e-15}  & \lower.3ex\hbox{18.57}  & \lower.3ex\hbox{2.692e-13} & \lower.3ex\hbox{12.90}  & \lower.3ex\hbox{3.476e-13}  & \lower.3ex\hbox{15.18} \\
512 & \lower.3ex\hbox{2.786e-15}  & \lower.3ex\hbox{-0.58}  & \lower.3ex\hbox{7.637e-13} & \lower.3ex\hbox{-1.50}  & \lower.3ex\hbox{1.046e-12}  & \lower.3ex\hbox{-1.59} \\
\hline
\end{tabular}
\end{center}
\label{SinTab}
\end{table}

\begin{table}
\caption{Relative errors and estimated orders of convergence for
interface I3}
\begin{center} \footnotesize
\begin{tabular}{ccccccc} \hline
& \multicolumn{2}{c}{Results for $\phi^n$} & \multicolumn{2}{c}{Results for $\partial\phi^n/\partial\tau$} & \multicolumn{2}{c}{Results for $\partial\phi_1^n/\partial\nu$}\\
\hline
$n$  &  Error & EOC &  Error & EOC & Error & EOC\\
\hline
16  & \lower.3ex\hbox{5.8841e-02}  & \lower.3ex\hbox{-}  & \lower.3ex\hbox{7.029e-02} & \lower.3ex\hbox{-}  & \lower.3ex\hbox{9.647e-01}     & \lower.3ex\hbox{-}  \\
32  & \lower.3ex\hbox{4.455e-03}  & \lower.3ex\hbox{3.72}  & \lower.3ex\hbox{7.084e-03} & \lower.3ex\hbox{3.31}  & \lower.3ex\hbox{6.920e-02}  & \lower.3ex\hbox{3.80}  \\
64  & \lower.3ex\hbox{2.517e-05}  & \lower.3ex\hbox{7.47}  & \lower.3ex\hbox{5.694e-05} & \lower.3ex\hbox{6.96}  & \lower.3ex\hbox{4.186e-04}  & \lower.3ex\hbox{7.37}  \\
128 & \lower.3ex\hbox{9.063e-10}  & \lower.3ex\hbox{14.76} & \lower.3ex\hbox{3.011e-09} & \lower.3ex\hbox{14.21}  & \lower.3ex\hbox{1.644e-08}  & \lower.3ex\hbox{14.63} \\
256 & \lower.3ex\hbox{1.847e-15}  & \lower.3ex\hbox{18.90}  & \lower.3ex\hbox{2.762e-13} & \lower.3ex\hbox{13.41}  & \lower.3ex\hbox{3.618e-13}  & \lower.3ex\hbox{15.47} \\
512 & \lower.3ex\hbox{2.715e-15}  & \lower.3ex\hbox{-0.56}  & \lower.3ex\hbox{8.407e-13} & \lower.3ex\hbox{-1.61}  & \lower.3ex\hbox{1.164e-12}  & \lower.3ex\hbox{-1.69} \\
\hline
\end{tabular}
\end{center}
\label{CosTab}
\end{table}

Tables \ref{SinTab} and \ref{CosTab} show the relative error and EOC
results for interfaces I2 and I3, respectively. Again
super-algebraic convergence is apparent until $n=256$, reaching
peaks of almost order 19 in both cases. Once again rounding errors
begin to destroy the convergence rate when $n=256$. Both higher
convergence rates and accuracy levels are observed in the
computations of $\phi^n$ as compared with it's directional
derivatives. However, since the derivatives are computed from
$\phi^n$ then this is probably to be expected; these trends are also
consistent with the observations made in \cite{MP11}.


\section{Conclusion}

The transmission problem for the Laplace equation on a periodic
half-space has been considered. The study was motivated by its
application to the modelling of electrified oil films used in the
development of novel switchable liquid optical devices (diffraction
gratings). A boundary integral formulation which reduces the problem
to the study of the interface alone was derived and solved in a
highly efficient manner using the Nystr\"{o}m method. The quadrature
rules were chosen with reference to supporting results in numerical
analysis, and were predicted to converge super-algebraically.
Numerical experiments demonstrated this convergence rate in practise
for a choice of parameters appropriate to our goal application, and
for a range of interface geometries.

\section*{Acknowledgments}
The author thanks Prof. Carl Brown, Dr. Reuben O'Dea and Dr. B.
Tomas Johansson for stimulating discussions. I am also indebted to
the anonymous reviewers for their careful reading, suggestions and
corrections.

\end{document}